\documentclass[preprint,12pt]{elsarticle}

\usepackage{amsmath,amsfonts,amssymb}
\usepackage{color}
\usepackage{graphicx}
\usepackage[utf8]{inputenc}
\usepackage[ruled,vlined]{algorithm2e}

\newproof{dem}{Proof}

\newtheorem{thm}{Theorem}[section]

\newcommand{\vect}[1]{\boldsymbol{#1}}

\journal{Journal of Computational and Applied Mathematics}

\begin{document}

\begin{frontmatter}



\title{An explicit-implicit Generalized Finite Difference scheme for a parabolic-elliptic density-suppressed motility system}

\cortext[cor1]{Corresponding author}

\author[ucm]{F. Herrero-Hervás}
\ead{fedher01@ucm.es}

\address[ucm]{Departamento de An\'alisis Matem\'atico y Matem\'{a}tica Aplicada, Instituto de Matemática Interdisciplinar, Universidad Complutense de Madrid, Madrid, Spain }
\begin{abstract}
In this work, a Generalized Finite Difference (GFD) scheme is presented for effectively computing the numerical solution of a parabolic-elliptic system modelling a bacterial strain with density-suppressed motility. The GFD method is a meshless method known for its simplicity for solving non-linear boundary value problems over irregular geometries. The paper first introduces the basic elements of the GFD method, and then an explicit-implicit scheme is derived. The convergence of the method is proven under a bound for the time step, and an algorithm is provided for its computational implementation. Finally, some examples are considered comparing the results obtained with a regular mesh and an irregular cloud of points.
\end{abstract}


\begin{keyword}
	Density-suppressed motility \sep Generalized finite difference method \sep Keller-Segel model
	
	
	
\end{keyword}

\end{frontmatter}

\section{Introduction}\label{intro}
Density-suppressed motility is a biological feature, introduced in 2011 in \cite{Liu}, through which the random motile motions of a strain of \textit{Escherichia coli} (\textit{E. coli}) cells are reduced at areas of high concentration of molecules of acyl-homoserine lactone (AHL). Moreover, the AHL is directly excreted by the bacteria and assumed to experience a time decayment. The process is modelled by the following non-linear system of parabolic equations over a domain $\Omega \subset \mathbb{R}^N$
\begin{equation} \label{1.1}
\begin{cases}
    \displaystyle \frac{\partial u}{\partial t} = \Delta (\gamma(v) u) + r u \left(1 - \frac{u}{K}\right), \\\\
    \displaystyle \tau \frac{\partial v}{\partial t} = D_v \Delta v - \alpha v + \beta u,
\end{cases} \quad \vect{x} \in \Omega, \hspace{0.1 cm} t >0,
\end{equation}
where $u$ represents the \textit{E. coli} cell density and $v$ the concentration of AHL. The function $\gamma: \mathbb{R} \to \mathbb{R}$ models the motility regulation, thus being monotone decreasing, though further assumptions are usually made to prove properties of the system. Parameters $r$ and $K$ encapsulate the logistic growth assumed for the bacteria, while $D_v$ is the diffusion coefficient of AHL and $\alpha$ and $\beta$ are respectively the decayment and production rates of AHL. Lastly $\tau \in \{0, 1\}$ distinguishes between fast diffusion scenarios (corresponding to $\tau = 0$) and regular diffusion.

When computing the Laplacian in the second equation of \eqref{1.1}, a special case of a Keller-Segel chemotaxis system (see \cite{KS1, KS2}) is retrieved, since $\Delta(\gamma(v) u) = \nabla \cdot \big [ \gamma(v) \nabla u + u \gamma'(v) \nabla v \big ]$. Here the diffusion and chemotaxis coefficients are $\gamma(v)$ and $\gamma'(v)$ respectively, hence being closely related.

Certain special properties arise when studying system \eqref{1.1}, such as a the formation of stripe patterns (experimentally obtained in \cite{Liu} and formally analyzed in \cite{Fu}) or the existence of travelling waves \cite{Li-TS}.

Our work is devoted to the numerical analysis of the parabolic-elliptic version of system \eqref{1.1}, this is, with $\tau = 0$ in a fast diffusion context. Nondimensionalizing the system we obtain
\begin{equation} \label{1.2}
\begin{cases}
    \displaystyle \frac{\partial u}{\partial t} = \Delta (\gamma(v) u) + \mu u (1 - u ), \\\\
    \displaystyle - \Delta v + v =  u,
\end{cases}  \quad \vect{x} \in \Omega, \hspace{0.1 cm} t >0.
\end{equation}
Analytically, this system has been studied in \cite{tello} under homogeneous Neumann boundary conditions. If the motility regulation function $\gamma$ verifies
\begin{equation} \label{1.3}
    \begin{cases}
        \gamma \in C^3([0, \infty)), \\\\
       \gamma(s) \geq 0, \hspace{0.2 cm} \gamma'(s) \leq 0, \hspace{0.2 cm} \gamma''(s) \geq 0, \hspace{0.2 cm} \gamma'''(s) \leq 0, \\\\
       -2 \gamma'(s) + \gamma''(s) s \leq \mu_0 < \mu,\\\\
       \displaystyle \frac{|\gamma'(s)|^2}{\gamma(s)} \leq c_\gamma \leq \infty,
    \end{cases} \hspace{0.5 cm} \text{for all } s\geq 0,
\end{equation}
and the initial value $u(\vect{x},0) = u_0(\vect{x})$ is such that
\begin{equation}\label{1.4}
\begin{cases}
    u_0 \in C^{2, \alpha}(\overline{\Omega}), \\\\
    \displaystyle \frac{\partial u_0}{\partial \Vec{n}} = 0 \text{ in } \partial \Omega \\\\
    0 < \underline{u}_0 < u_0 < \overline{u}_0 < \infty,
\end{cases}    
\end{equation}
for certain positive constants $\mu_0$, $c_\gamma$, $\overline{u}_0$, $\underline{u}_0$, then it is proven that a unique solution globally exists in time, satisfying:
\begin{equation} \label{1.5}
\lim_{t\to\infty} ||u-1||_{L^\infty(\Omega)} + ||v-1||_{L^\infty(\Omega)} = 0.    
\end{equation}
In this paper we study a Generalized Finite Difference (GFD) Method to numerically solve system \eqref{1.2} together with Neumann homogeneous boundary conditions. The GFD method was initially proposed in 1960 by Collatz in \cite{collatz} and by Forsythe and Wasow in \cite{forsythe}, though it was later reintroduced in the 1970s by Jensen \cite{jensen} and Perrone and Kao \cite{perrone}. Over the years many other authors have studied improvements of the method, as seen for example in \cite{Benito-factors}.

The GDF method consists of obtaining finite-difference approximation formulae for the partial derivatives, based however on irregular clouds of points. Thus, the method allows for geometrically irregular domains, where the classical Finite Difference Method cannot be used, as no regular grid can be considered.

Several recent papers have considered the GFD solution to nonlinear systems arising in physical and biological processes, such as \cite{Benito-GFD}, where the convergence to periodic solutions of a chemotaxis system is considered; \cite{Benito-tumor}, taking into account a model for tumor growth involving nutrient density, extracellular matrix and matrix degrading enzymes; or \cite{Waterwaves}, where several types of water waves are simulated.

The article is structured as follows: firstly, for it to be self-contained, Section \ref{s1} begins with an introduction of the GFD method, followed by the proposed numerical scheme, an algorithm for its computational implementation and the main result of the paper, Theorem \ref{num1}, regarding the convergence of the method. Then, in Section \ref{s2} numerical tests are performed to assess the results of the method. Lastly, the conclusions are outlined in Section \ref{s3}.

\section{Explicit-implicit GFD Scheme} \label{s1}
\subsection{Brief description of the GFD Method} \label{preliminares}
To account for general systems, we first begin by considering the general setting for the GFD Method, later allowing us to deduce an explicit-implicit scheme for numerically solving \eqref{1.2}. To do so, we consider a problem of the form:
\begin{equation} \label{2.1}
\begin{cases}
\displaystyle  \frac{\partial u}{\partial t} (\vect{x},t) = L_\Omega [u(\vect{x},t)], \quad \text{for } \vect{x} \in \Omega, \hspace{0.1 cm} t>0, \\\\
\displaystyle  \frac{\partial u}{\partial \vec{n}}  (\vect{x},t) = 0, \quad \text{for } \vect{x} \in \partial \Omega, \hspace{0.1 cm}  t>0, \\\\
\displaystyle  u(\vect{x},0) = g(\vect{x}), \quad \text{for } \vect{x} \in \Omega,
\end{cases}
\end{equation}
over a bounded domain $\Omega \subset \mathbb{R}^N$, where $L_\Omega$ is a second order differential operator defined in $\Omega$, this is, only including $\partial / \partial x$, $\partial / \partial y$, $\partial^2 / \partial x^2$, $\partial^2 / \partial y^2$ and $\partial^2/\partial x \partial y$.

The strategy of the method is to consider a time discretization of constant step $\Delta t >0$ with which we can consider a first order approximation of the time derivative $\frac{\partial}{\partial t}$. The use of a forward difference formula yields an explicit scheme. The spatial derivatives are approximated using finite difference formulae that have to be derived.

More precisely, we consider a discretization of $\Omega$ made up of $m$ \textit{nodes}, $M = \{z_1, \dots, z_m\} \subset \Omega$ over which we seek to approximate the solution to the problem. These are the points that will be used to obtain the finite difference formulae for the spatial derivatives. As usual, we denote by $U_j^n$ the approximation of $u(\vect{z_j}, n \Delta t)$. To keep the notation simple we assume that the domain $\Omega$ is of two spatial dimensions, though higher dimensions are analogously treated. This way, we consider both coordinates of each node given by $\vect{z_j} = (x_j, y_j)$.

\begin{figure}[htp]
	\centering
	\includegraphics[scale=0.35]{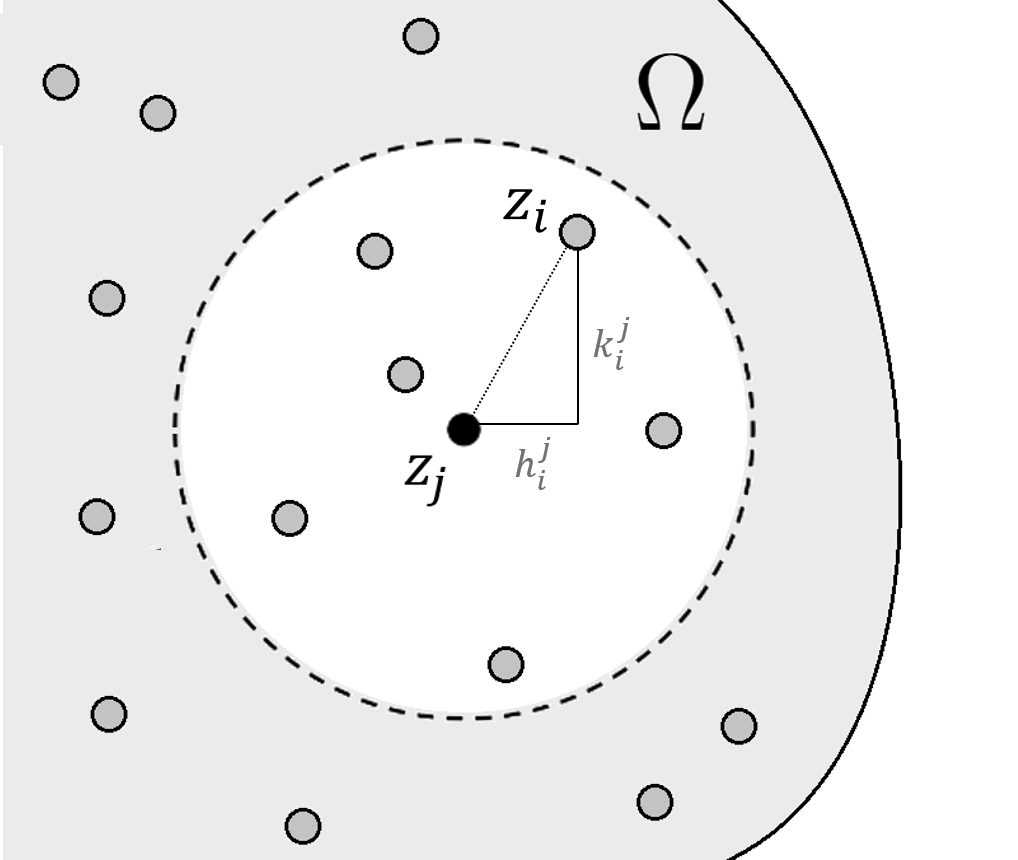}
	\caption{Representation of the nodes that make up an $E_6-$star centered in $\vect{z_j}$.}
	\label{fig:f1}
\end{figure}

For every node $\vect{z_j}$ and for a fixed integer $s$, we define an $E_s-$star centered in $\vect{z_j}$ as a set of $s$ other nodes in $M$ located within a neighbourhood of $\vect{z_j}$. For all the nodes in this star, $\vect{z_i} = (x_i, y_i) \in E_s$, we define 
$$h_i^j = x_i - x_j, \quad k_i^j = y_i - y_j.$$
Figure \ref{fig:f1} depicts an $E_6-$star centered in $\vect{z_j}$ constituted by the nodes inside the white neighbourhood of $\vect{z_j}$. For the node $\vect{z_i}$, the values of $h_i^j$ and $k_i^j$ are also represented.
As the main goal is to obtain approximation formulae for the spatial  derivatives, a Taylor approximation is considered. Since the operator $L_\Omega$ was assumed to be of second order, we only require a second order Taylor expansion. This way, we have:
\begin{equation} \label{2.2}
\begin{split}
&u(x_i, y_i) \approx u(x_j, y_j) + h_i^j \frac{\partial u}{\partial x}(x_j, y_j) + k_i^j \frac{\partial u}{\partial y}(x_j, y_j) \\&+ \frac{1}{2}\left[(h_i^j)^2 \frac{\partial^2 u}{\partial x^2}(x_j, y_j) +  (k_i^j)^2 \frac{\partial^2 u}{\partial y^2}(x_j, y_j) + 2h_i^j k_i^j \frac{\partial^2 u}{\partial x \partial y}(x_j, y_j)  \right].
\end{split}
\end{equation} 
And thus, to obtain the best possible formulae for the partial derivatives, the following function is considered, being a weighted sum of squares of the differences between the value of every $U_i^n$ and its previously obtained second order Taylor approximation centered in $\vect{z_j}$:
\begin{equation} \label{2.3}
\begin{split}
B(U_j^n) =& \sum_{i=1}^s (w_i^j)^2 \Biggl[ U_j^n + h_i^j \frac{\partial U_j^n}{\partial x} +  k_i^j \frac{\partial U_j^n}{\partial y} \\&+ \frac{1}{2} \left( (h_i^j)^2 \frac{\partial^2  U_j^n}{\partial x^2} +  (k_i^j)^2 \frac{\partial^2  U_j^n}{\partial y^2} + 2h_i^j k_i^j \frac{\partial^2  U_j^n}{\partial x \partial y} \right) - U_i^n  \Biggr]^2.
\end{split}
\end{equation}
The partial derivatives of $U_j^n$ denote their approximations, which we seek to determine. Moreover $w_i^j$ are the weights, non-negative symmetric functions that decrease with the distance from $\vect{z_j}$ to $\vect{z_i}$, usually given by:
$$w_i^j = \frac{1}{(h_i^j + k_i^j)^{\alpha/2}} = \frac{1}{||\vect{z_j} - \vect{z_i}||^\alpha},$$
for a certain $\alpha > 0$. 

This way, minimizing the sum $B$ with respect to the partial derivatives of $U_j^n$ yields the best finite difference formulae for them. Moreover, $B$ being a sum of quadratic terms leads us to obtain its minimum through the solution of a linear system of the form:
\begin{equation} \label{2.4}
\vect{A_j^n} \cdot \vect{D_j^n} = \vect{b_j^n},
\end{equation}
where:
\begin{equation} \label{2.5}
\vect{A_j^n} = \left(
\begin{array}{cccc}
h^j_1 & h^j_2 & \cdots & h^j_s \\
k^j_1 & k^j_2 & \cdots & k^j_s \\
\vdots & \vdots & \vdots & \vdots \\
h^j_1k^j_1 & h^j_2k^j_2 & \cdots & h^j_sk^j_s \\
\end{array}
\right)\left(
\begin{array}{cccc}
(w_1^j)^2 &  &  &  \\
& (w_2^j)^2 &  &  \\
&  & \cdots &  \\
&  &  & (w_s^j)^2 \\
\end{array}
\right)\left(
\begin{array}{cccc}
h^j_1 & k^j_1 & \cdots & h^j_1k^j_1 \\
h^j_2 & k^j_2 & \cdots & h^j_2k^j_2 \\
\vdots & \vdots& \vdots & \vdots \\
h^j_s & k^j_s &  \cdots & h^j_sk^j_s \\
\end{array}\right)  
\end{equation}
\vspace{0.02 cm}
\begin{equation} \label{2.6}
%\small
\vect{b_j^n} = \begin{pmatrix}
\displaystyle    \sum_{i=1}^s (U_i^n - U_j^n)h_i^j (w_i^j)^2 \\\\ \displaystyle \sum_{i=1}^s (U_i^n - U_j^n)k_i^j (w_i^j)^2 \\\\ \displaystyle \frac{1}{2} \sum_{i=1}^s (U_i^n - U_j^n)(h_i^j)^2 (w_i^j)^2 \\\\ \displaystyle
\frac{1}{2} \sum_{i=1}^s (U_i^n - U_j^n)(k_i^j)^2 (w_i^j)^2 \\\\ \displaystyle
\sum_{i=1}^s (U_i^n - U_j^n)h_i^j k_i^j (w_i^j)^2
\end{pmatrix},
\end{equation}
and $\vect{D_j^n}$ is the unknown vector, containing the partial derivatives:
\begin{equation} \label{2.7}
\vect{D_j^n} = \begin{pmatrix} \displaystyle 
\frac{\partial U_j^n}{\partial x} & \displaystyle \frac{\partial U_j^n}{\partial y} &  \displaystyle \frac{\partial^2 U_j^n}{\partial x^2} &  \displaystyle \frac{\partial^2 U_j^n}{\partial y^2} &  \displaystyle \frac{\partial^2 U_j^n}{\partial x \partial y}
\end{pmatrix}^T.
\end{equation}
Consequently, the solution to system \eqref{2.4} provides the sought finite difference formulae for $\partial U_j^n / \partial x$, $\partial U_j^n/\partial y$, $\partial^2 U_j^n/\partial x^2$, $\partial^2 U_j^n / \partial y^2$ and $\partial^2 U_j^n / \partial x \partial y$, obtained through the values of $u$ on the nodes of the star.

Regarding the time derivative, the above mentioned forward difference scheme yields:
\begin{equation} \label{2.8}
    \frac{\partial u}{\partial t}(x_j, y_j, n \Delta t) \approx \frac{U_j^{n+1} - U_j^n}{\Delta t}.
\end{equation}
Therefore, to obtain the GFD scheme, the partial derivatives in $L_\Omega$ in \eqref{2.1} are replaced by their finite difference approximations given as the solution to system \eqref{2.4} for the spatial derivatives and as in \eqref{2.8} for the time derivative.

To obtain the GFD scheme, it is however desirable to express the solution to system \eqref{2.4} directly as sum of the values of $u$ over the nodes, just as in the classical Finite Difference Method. To simplify the notation with the sub and super indexes, instead of a general node $\vect{z_j} = (x_j, y_j)$ we consider a fixed node $\vect{z_0} = (x_0, z_0)$, denoting $\vect{A_j^n}$ just by $\vect{A}$ and the distances $h_i^j$ and $k_i^j$ and weights $w_i^j$ by $h_i$, $k_i$ and $w_i$. For the desired representation formula, for each node in the $E_s-$star centered in $\vect{z_0}$ we consider the vector:
\begin{equation} \label{2.9}
\vect{c_i} = \begin{pmatrix}
h_i & k_i & \frac{1}{2} h_i^2 &\frac{1}{2} k_i^2 &  h_i k_i
\end{pmatrix}.
\end{equation}
Hence, it follows that the solution to system \eqref{2.4} can be expressed as:
\begin{equation} \label{2.10}
\begin{cases}
\displaystyle   \frac{\partial U(x_0,y_0, n \Delta t)}{\partial x} = - \lambda_{01} U_0^n + \sum_{i=1}^s \lambda_{i1}U_i^n , \\\\\
\displaystyle \frac{\partial U(x_0,y_0, n \Delta t)}{\partial y} = - \lambda_{02} U_0^n + \sum_{i=1}^s \lambda_{i2}U_i^n ,\\\\
\displaystyle \frac{\partial^2 U(x_0,y_0, n \Delta t)}{\partial x^2} = - \lambda_{03} U_0^n + \sum_{i=1}^s \lambda_{i3}U_i^n ,\\\\
\displaystyle \frac{\partial^2 U(x_0,y_0, n \Delta t)}{\partial y^2} = - \lambda_{04} U_0^n + \sum_{i=1}^s \lambda_{i4}U_i^n , \\\\
\displaystyle \frac{\partial^2 U(x_0,y_0, n \Delta t)}{\partial x \partial y} = - \lambda_{05} U_0^n + \sum_{i=1}^s \lambda_{i5}U_i^n ,
\end{cases}
\end{equation}
where the finite difference coefficients $\lambda_{i,j}$ are given by:
\begin{equation} \label{2.11}
\lambda_{ir} = w_i^2 (\vect{A}^{-1} \vect{c_i})_r, \quad \lambda_{0i} = \sum_{i=1}^s \lambda_{i},
\end{equation}
for $r\in \{1, \dots, s\}$, $i \in \{1, \dots, 5\}$, being $(\vect{A}^{-1} \vect{c_i})_r$ the $r-$th coordinate of the vector $\vect{A}^{-1} \vect{c_i}$. Given that the Laplacian will often appear in the scheme, we denote $\lambda_{00} = \lambda_{03} + \lambda_{04}$ and $\lambda_{i0} = \lambda_{i3} + \lambda_{i4}$.

Direct substitution of the finite difference formulae for the derivatives in \eqref{2.1} leads to the corresponding explicit GFD scheme, of first order in time and second order in space. As a remark, $\vect{A}^{-1}$ can be computed through a Cholesky factorization, given that $\vect{A}$ is a positive definite matrix (see, for example \cite{R6,vargas}).

\begin{figure}[htp]
    \centering
    \includegraphics[scale=0.35]{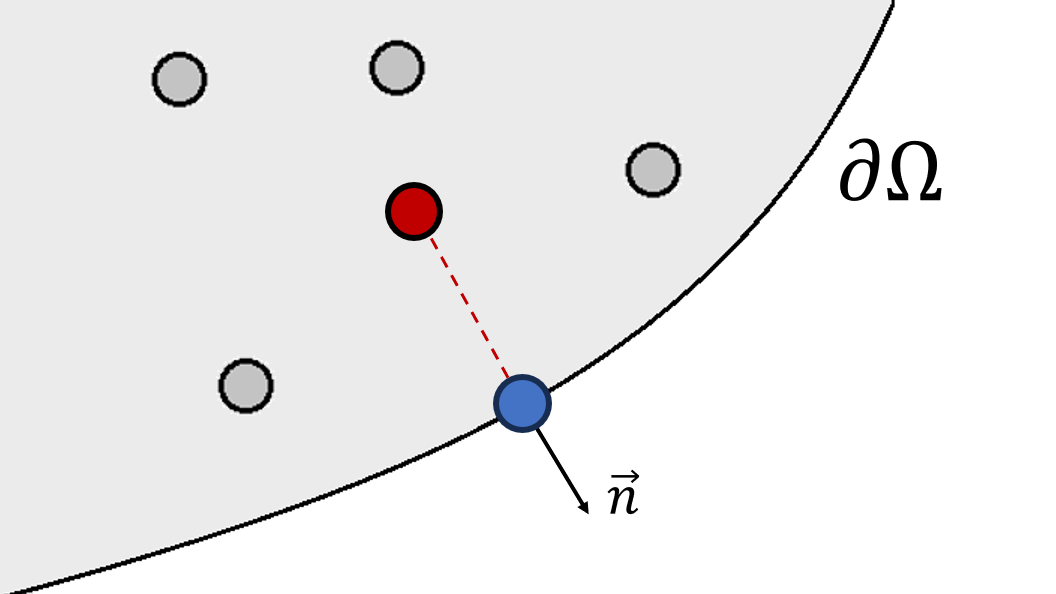}
    \caption{A boundary node (colored in blue) and an inner node (in red) placed along the direction of the normal vector.}
    \label{fig:f2}
\end{figure}

The boundary conditions are equally treated. In this case, for the Neumann condition on the boundary nodes the normal derivative is approximated through its surrounding nodes or directly with a first order scheme, placing an inner node in the direction of the normal vector, as schematically represented in Figure \ref{fig:f2}. The inner node coloured in red can be used to obtain a first order approximation of the normal derivative on the blue boundary node.
\subsection{Scheme and algorithm for the density-suppressed motility model}
For system \eqref{1.2}, the procedure described above can be easily adapted to account for the second equation with no time derivative. Direct substitution of the finite difference formulae for the derivatives leads to the following scheme:
\begin{equation}
\label{2.12}
\left\lbrace
\begin{aligned}{}
&\frac{U_0^{n+1}-U_0^n}{\Delta t}=\gamma(V^{n}_0)\left[-\lambda_{00}U_0^n+\sum^{s}_{i=1}\lambda_{i0}U_i^n\right] \\&\qquad\qquad
+2\gamma'(V^n_0)\left(-\lambda_{01}U_0^n+\sum^{s}_{i=1}\lambda_{i1}U_i^n\right)\left(-\lambda_{01}V_0^n+\sum^{s}_{i=1}\lambda_{i1}V_i^n\right)
\\&\qquad\qquad  +2\gamma'(V^n_0)\left(-\lambda_{02}U_0^n+\sum^{s}_{i=1}\lambda_{i2}U_i^n\right)\left(-\lambda_{02}V_0^n+\sum^{s}_{i=1}\lambda_{i2}V_i^n\right) 
\\&\qquad\qquad+ U^n_0\gamma''(V^n_0)\left[\left(-\lambda_{01}V_0^n+\sum^{s}_{i=1}\lambda_{i1}V_i^n\right)^2+\left(-\lambda_{02}V_0^n+\sum^{s}_{i=1}\lambda_{i2}V_i^n\right)^2\hspace{0.05 cm}\right]\\
&\qquad\qquad+ U^n_0\gamma'(V^n_0)\Big( V^n_0 -U^n_0 \Big) + \mu U^n_0(1-U^n_0) + \mathcal{O}(\Delta t,h_i^2,k_i^2),
\end{aligned}
\right.
\end{equation} and
\begin{equation}\label{2.13}
V_0^n-\left[-\lambda_{00}V_0^n+\sum^{s}_{i=1}\lambda_{i0}V_i^n\right]=U^n_0 +\mathcal{O}(\Delta t,h_i^2,k_i^2).
\end{equation}
The resulting scheme is explicit-implicit, as it is initialized in $U_0^0$ with the values of $u_0(\vect{x})$, which are needed to compute $V_0^0$ by solving the resulting system \eqref{2.13}. These values are used for obtaining first $U_0^1$, explicitly through \eqref{2.12} and then secondly $V_0^1$, again implicitly. The process is repeated until a chosen final time.

To give a more detailed description of the method, an algorithm for the computations is provided below.
\\
\vspace{0.5 cm}
\\
\begin{algorithm}[H]
\SetAlgoLined
\vspace{0.2 cm}
\underline{\textbf{Initialization}}\\
Import the set of nodes (\texttt{NT} $\leftarrow$ total number of nodes)\\
Select number of nodes per star, \texttt{n\_nodes}; time step and weights, $w$\\
Create empty vectors to store the values of $u$ and $v$\\
$\rightarrow$ The vector for $u$ is initialized with its initial condition $u(\vect{x},0)$ \\
\vspace{0.2 cm}
\underline{\textbf{Star generation}} \\
Distinguish between inner and boundary nodes
\\ \For{i = 1:\texttt{NT}}{
\If{node $i$ is an inner node}{
\For{j = 1:\texttt{NT}}{
Calculate distance between both nodes
}
Select the \texttt{n\_nodes} closest nodes\\
\% This can be done by storing a matrix with $1$ on position $(i,j)$ if node $j$ is on the star centered in node $i$ and $0$ otherwise 
}
}
\vspace{0.2 cm}

\end{algorithm}

\begin{algorithm}[H]
\SetAlgoLined
\vspace{0.2 cm}
\underline{\textbf{Finite difference approximations of the derivatives}} \\
\texttt{NI} $\leftarrow$ total number of inner nodes\\
$D1 = D2 = D3 = D4 = D5 = zeros(\texttt{NI}, \texttt{n\_nodes})$ 
\\ \% These are the matrices where the coefficients for the partial derivatives will be stored \\
\For{i = 1: \texttt{NI}}{
A = zeros(5,5) \\
counter = 0 \\
\For{j = 1:\texttt{NT}}{
\If{Node $j$ is in the star centered in node $i$}{
Compute $h_i^j$, $k_i^j$ and $w_i^j$ \\
$c = \left(h_i^j, k_i^j, (h_i^j)^2/2, (k_i^j)^2/2, h_i^jk_i^j/2 \right)$ \\
$b_{aux}(:, counter) = c^T\cdot w_i^j$ \\
$A = A + c^T \cdot c \cdot w_i^j$\\
$counter = counter + 1$
}
} \vspace{0.2 cm}
\% We obtain $D_{aux}$ through the Cholesky factorization of $A$ \\
$ R \leftarrow$ Cholesky decomposition of $A$ (lower triangular matrix) \\
$ M = R \setminus eye(5)$ \\
$ D_{aux} = M^T \cdot M + b_{aux}$
\\
\vspace{0.2 cm} \% We store the partial derivative coefficients for each node \\
$D1(i,:) \leftarrow $ 1st row of $D_{aux}$ \% coefficients for $\partial / \partial x$ \\
$D2(i,:) \leftarrow $ 2nd row of $D_{aux}$ \% coefficients for $\partial / \partial y$ \\
$D3(i,:) \leftarrow $ 3rd row of $D_{aux}$ \% coefficients for $\partial^2 / \partial x^2$ \\
$D4(i,:) \leftarrow $ 4th row of $D_{aux}$ \% coefficients for $\partial^2 / \partial y^2$ \\
$D5(i,:) \leftarrow $ 5th row of $D_{aux}$ \% coefficients for $\partial^2 / \partial x \partial y$ \\
}
\end{algorithm}

\begin{algorithm}[H]
\SetAlgoLined
\vspace{0.2 cm}
\underline{\textbf{Temporal loop}} \\
$J \leftarrow$ matrix with the coefficients needed to solve the elliptic equation on each time step \\
$sol\_u \leftarrow$ vector where the solution $u$ will be stored on each time step. Initialized with its initial condition $u(x,0)$ \\
$sol\_v \leftarrow$ empty vector where the solution $v$ will be stored on each time step. Does not require an initial value \\
\texttt{step\_num} $\leftarrow$ number of time steps obtained from the selected $\Delta t$ and final time\\
\For{n = 1:\texttt{step\_num}}{
$sol\_v \leftarrow$ solve the linear system for the elliptic equation \\
\% For the parabolic equation we use the explicit scheme presented in the paper \\
\For{i = 1: \texttt{NI}}{
\For{j = 1:\texttt{n\_nodes}}{
Compute the differences $u_{dif}(j) = U_i - U_j$, $v_{dif}(j) = V_i - V_j$ with the current values of $sol\_u$ and $sol\_v$ \\
\% These are needed to obtain the finite difference approximation formulae through the already calculated coefficients stored in $D1 \dots D5$
}
$new\_sol\_u \leftarrow$ $U_i^{n+1}$ computed via the explicit scheme through the finite difference approximations of the partial derivatives based on the values of $U_i^n$ and $V_i^n$
}
\% Boundary nodes \\
For the Neumann boundary conditions, a simple forward or backward difference scheme can be used \\
\vspace{0.2 cm}
$sol\_u = new\_sol\_u$ \% Update the solution for the next time step
}
    
\end{algorithm}

\subsection{Convergence of the method}
As per usual with explicit schemes, certain condition shall be satisfied to guarantee its convergence. The main result is established in Theorem \ref{num1}.
\begin{thm} \label{num1}
	Let $(u,v)$ be the solution to system \eqref{1.2} under conditions \eqref{1.3}-\eqref{1.4}, then the GFD scheme \eqref{2.9}-\eqref{2.10} is convergent if the time step $\Delta t$ is such that
	$$\Delta t<\dfrac{1+|1-\lambda_{00}|+\sum_{i=1}^s|\lambda_{i0}|}{\Biggl[|1-\lambda_{00}|+\sum_{i=1}^s|\lambda_{i0}|\Biggr](A'_1+A''_1)+B_1},$$
	holds for every inner node, where coefficients $A'_1$, $A''_1$ and $B_1$ are stated in the proof.
\end{thm}
\begin{dem}
To establish conditions for the convergence, we consider $eu^n_0:= U^n_0-u^n_0$ the difference between the discrete and the continuous solution at the node $\vect{z_0}$, defining $eu^n_i,$ and $ev^n_i$ analogously. Taking into account that the exact values of $u$ and $v$ must satisfy system \eqref{1.2}, the main proof strategy is to obtain an expression that $eu^n_0$ verifies and conditions for its convergence to 0.

For example, for the fist term of the scheme, by the Mean Value Theorem, we have that:
	\begin{equation}\label{prueba}
	\begin{split}
	\gamma&(V^n_0)\left[-\lambda_{00}U_0^n+\sum^{s}_{i=1}\lambda_{i0}U_i^n\right]-\gamma(v^n_0)\left[-\lambda_{00}u_0^n+\sum^{s}_{i=1}\lambda_{i0}u_i^n\right] \\&\pm\gamma(V^n_0)\left[-\lambda_{00}u_0^n+\sum^{s}_{i=1}\lambda_{i0}u_i^n\right]=\\
	&=\gamma'(\xi)\left[-\lambda_{00}u_0^n+\sum^{s}_{i=1}\lambda_{i0}u_i^n\right]ev^n_0+\gamma(V^n_0)\left[-\lambda_{00}eu_0^n+\sum^{s}_{i=1}\lambda_{i0}eu_i^n\right],
	\end{split}
	\end{equation}
After performing similarly in the rest of the terms, it is obtained that:
	\begin{equation}\label{prueba7}
	\begin{split}
	\small &\frac{eu^{n+1}_0-eu^n_0}{\Delta t}=\gamma''(\xi)\left[-\lambda_{00}u_0^n+\sum^{s}_{i=1}\lambda_{i0}u_i^n\right]ev^n_0+\gamma'(V^n_0)\left[-\lambda_{00}eu_0^n+\sum^{s}_{i=1}\lambda_{i0}eu_i^n\right]\\
	&+2\gamma''(\xi)\left(-\lambda_{01}U_0^n+\sum^{s}_{i=1}\lambda_{i1}U_i^n\right)\left(-\lambda_{01}V_0^n+\sum^{s}_{i=1}\lambda_{i1}V_i^n\right)ev^n_0+\\
	&+2\gamma'(V^n_0)\left(-\lambda_{01}u_0^n+\sum^{s}_{i=1}\lambda_{i1}u_i^n\right)\left(-\lambda_{01}ev_0^n+\sum^{s}_{i=1}\lambda_{i1}ev_i^n\right)+\\
	&+2\gamma'(V^n_0)\left(-\lambda_{01}eu_0^n+\sum^{s}_{i=1}\lambda_{i1}eu_i^n\right)\left(-\lambda_{01}V_0^n+\sum^{s}_{i=1}\lambda_{i1}V_i^n\right)+\\
	&+2\gamma''(\xi)\left(-\lambda_{02}U_0^n+\sum^{s}_{i=1}\lambda_{i2}U_i^n\right)\left(-\lambda_{02}V_0^n+\sum^{s}_{i=1}\lambda_{i2}V_i^n\right)ev^n_0+\\
	&+2\gamma'(V^n_0)\left(-\lambda_{02}eu_0^n+\sum^{s}_{i=1}\lambda_{i2}eu_i^n\right)\left(-\lambda_{02}V_0^n+\sum^{s}_{i=1}\lambda_{i2}V_i^n\right)+\\
	&+2\gamma'(V^n_0)\left(-\lambda_{02}eu_0^n+\sum^{s}_{i=1}\lambda_{i2}eu_i^n\right)\left(-\lambda_{02}V_0^n+\sum^{s}_{i=1}\lambda_{i2}V_i^n\right)+\\
	&+eu^n_0\gamma''(V^n_0)\Biggl(-\lambda_{01}V_0^n+\sum^{s}_{i=1}\lambda_{i1}V_i^n\Biggr)^2+eu^n_0\gamma''(V^n_0)\Biggl(-\lambda_{02}V_0^n+\sum^{s}_{i=1}\lambda_{i2}V_i^n\Biggr)^2\\
	&u^n_0\gamma'''(\xi)\Biggl(-\lambda_{01}V_0^n+\sum^{s}_{i=1}\lambda_{i1}V_i^n\Biggr)^2ev^n_0+u^n_0\gamma'''(\xi)\Biggl(-\lambda_{02}V_0^n+\sum^{s}_{i=1}\lambda_{i2}V_i^n\Biggr)^2ev^n_0\\
	&+u^n_0\gamma''(V^n_0)\Biggl(-\lambda_{01}ev_0^n+\sum^{s}_{i=1}\lambda_{i1}ev_i^n\Biggr)\Biggl(-\lambda_{01}(v_0^n+V^n_0)+\sum^{s}_{i=1}\lambda_{i1}(v_i^n+V^n_i)\Biggr)\\
	&+u^n_0\gamma''(V^n_0)\Biggl(-\lambda_{02}ev_0^n+\sum^{s}_{i=1}\lambda_{i2}ev_i^n\Biggr)\Biggl(-\lambda_{02}(v_0^n+V^n_0)+\sum^{s}_{i=1}\lambda_{i2}(v_i^n+V^n_i)\Biggr)\\
	&+eu^n_0V^n_0\gamma'(V^n_0)+u^n_0v^n_0\gamma''(\xi)ev^n_0+u^n_0ev^n_0\gamma'(V^n_0)-(U^n_0)^2\gamma''(\xi)ev^n_0\\
	&-eu^n_0(u^n_0+U^n_0)\gamma'(v^n_0)+\mu eu^n_0-\mu eu^n_0(u^n_0+U^n_0)
	\end{split}
	\end{equation}
	If we take $eu^n=\displaystyle\max_{i=0,...,s}\{|eu^n_i|\}$ and $ev^n = \displaystyle\max_{i=0,...,s}\{|ev^n_i|\}$ a bound for the previous expression can be obtained, of the form:
	\begin{equation}\label{prueba8}
	eu^{n+1}\leq A_1eu^n+B_1ev^n,
	\end{equation}
	with positive coefficients $A_1,$ and $ B_1$ that depend on $\Delta t$:
	\begin{equation}\label{prueba9}
	\begin{split}
	A_1:&=\Biggl|1+ \Delta t\Biggl[-\lambda_{00}-2\gamma'(V^n_0)\lambda_{01}\Biggl(-\lambda_{01}V^n_0+\sum_{i=1}^s\lambda_{i1}V^n_i\Biggr)\\
	&-2\gamma'(V^n_0)\lambda_{02}\Biggl(-\lambda_{02}V^n_0+\sum_{i=1}^s\lambda_{i2}V^n_i\Biggr)+\gamma''(V^n_0)\Biggl(-\lambda_{01}V^n_0+\sum_{i=1}^s\lambda_{i1}V^n_i\Biggr)^2\\
	&+\gamma''(V^n_0)\Biggl(-\lambda_{02}V^n_0+\sum_{i=1}^s\lambda_{i2}V^n_i\Biggr)^2+V^n_0\gamma'(V^n_0)+u^n_0\gamma'(V^n_0)\\
	&-(u^n_0+U^n_0)\gamma'(v^n_0)+\mu-\mu(u^n_0+U^n_0)\Biggr]\Biggr|+\Delta t\Biggl[|\gamma'(v^n_0)\sum_{i=1}^s\lambda_{i0}|\\
	&+2|\gamma''(v^n_0)\Biggl(-\lambda_{01}V^n_0+\sum_{i=1}^s\lambda_{i1}V^n_i\Biggr)\sum_{i=1}^s\lambda_{i1}|\\&+2|\gamma''(v^n_0)\Biggl(-\lambda_{02}V^n_0+\sum_{i=1}^s\lambda_{i2}V^n_i\Biggr)\sum_{i=1}^s\lambda_{i2}|\Biggr],
	\end{split}
	\end{equation}
	\vspace{0.3 cm}
	\begin{equation}\label{prueba10}
	\begin{split}
	B_1:&=\Delta t\Biggl|\gamma'(\xi)\Biggl(-\lambda_{00}U^n_0+\sum_{i=1}^s\lambda_{i0}U^n_i\Biggr)\\&+2\gamma''(\xi)\Biggl(-\lambda_{01}U^n_0+\sum_{i=1}^s\lambda_{i1}U^n_i\Biggr)\Biggl(-\lambda_{01}V^n_0+\sum_{i=1}^s\lambda_{i1}V^n_i\Biggr)\\
	&+2\gamma''(\xi)\Biggl(-\lambda_{02}U^n_0+\sum_{i=1}^s\lambda_{i2}U^n_i\Biggr)\Biggl(-\lambda_{02}V^n_0+\sum_{i=1}^s\lambda_{i2}V^n_i\Biggr)\\&-2\gamma'(V^n_0)\Biggl(-\lambda_{01}U^n_0+\sum_{i=1}^s\lambda_{i1}U^n_i\Biggr)\lambda_{01}
	-2\gamma'(V^n_0)\Biggl(-\lambda_{02}U^n_0+\sum_{i=1}^s\lambda_{i2}U^n_i\Biggr)\lambda_{02}\\&+u^n_0\gamma'''(\xi)\Biggl(-\lambda_{01}V^n_0+\sum_{i=1}^s\lambda_{i1}V^n_i\Biggr)^2 +u^n_0\gamma'''(\xi)\Biggl(-\lambda_{02}V^n_0+\sum_{i=1}^s\lambda_{i2}V^n_i\Biggr)^2 \\&-u^n_0\gamma''(v^n_0)\Biggl(-\lambda_{01}(v^n_0+V^n_0)+\sum_{i=1}^s\lambda_{i1}(v^n_i+V^n_i)\Biggr)\lambda_{01}\\
	&-u^n_0\gamma''(v^n_0)\Biggl(-\lambda_{02}(v^n_0+V^n_0)+\sum_{i=1}^s\lambda_{i2}(v^n_i+V^n_i)\Biggr)\lambda_{02}\\&+u^n_0v^n_0\gamma''(\xi)-(U^n_0)^2\gamma''(\xi)\Biggr|+\Delta t\Biggl[2|\gamma'(v^n_0)\Biggl(-\lambda_{01}U^n_0+\sum_{i=1}^s\lambda_{i1}U^n_i\Biggr)\sum_{i=1}^s\lambda_{i1}|\\&+2|\gamma'(v^n_0)\Biggl(-\lambda_{02}U^n_0+\sum_{i=1}^s\lambda_{i2}U^n_i\Biggr)\sum_{i=1}^s\lambda_{i2}|\\&+|u^n_0\gamma''(v^n_0)\Biggl(-\lambda_{01}(v^n_0+V^n_0)+\sum_{i=1}^s\lambda_{i1}(v^n_i+V^n_i)\Biggr)\sum_{i=1}^s\lambda_{i1}|\\&+|u^n_0\gamma''(v^n_0)\Biggl(-\lambda_{02}(v^n_0+V^n_0)+\sum_{i=1}^s\lambda_{i2}(v^n_i+V^n_i)\Biggr)\sum_{i=1}^s\lambda_{i2}| \Biggr].
	\end{split}
	\end{equation}
	Moreover, $A_1$ can be rewritten as $A_1=|1-\Delta t A_1^{'}|+\Delta t A_{1}^{''}$, for a clear choice of $A_1^{'}$ and $A_1^{''}$.
	\\\\Following a similar procedure in the scheme for $v$, we obtain:
	\begin{equation}\label{prueba13}
	eu^{n} \leq \Biggl[|1-\lambda_{00}|+\sum_{i=1}^s|\lambda_{i0}|\Biggr] ev^n.
	\end{equation}
	Substitution of \eqref{prueba13} in \eqref{prueba8} for $eu^{n}$ yields
	\begin{equation}\label{prueba13b}
	eu^{n+1}\leq \Biggl(A_1\Biggl[|1-\lambda_{00}|+\sum_{i=1}^s|\lambda_{i0}|\Biggr]+B_1\Biggr)ev^n.
	\end{equation}
	 For the desired convergence, we have to ensure that $eu^{n}$ and $ev^{n}$ tend to 0, therefore we impose
	\begin{equation} \label{prueba14}
	A_1[|1-\lambda_{00}|+\sum_{i=1}^s|\lambda_{i0}|]+B_1 < 1. 
	\end{equation}
	And due to the fact that $A_1=|1-\Delta t A_1^{'}|+\Delta t A_{1}^{''}$, we obtain the following bound for $\Delta t$:
	\begin{equation} \label{prueba15}
	\Delta t<\dfrac{1+|1 - \lambda_{00}|+\sum_{i=1}^s|\lambda_{i0}|}{\Biggl[|1-\lambda_{00}|+\sum_{i=1}^s|\lambda_{i0}|\Biggr](A'_1+A''_1)+B_1}.
	\end{equation}
\end{dem}
\section{Numerical Tests} \label{s2}
Lastly, we devote this section to show numerical solutions of system \eqref{1.2}, computed through the GFD scheme \eqref{2.12}-\eqref{2.13}, with different parameter values and initial conditions. To check the asymptotic convergence to the state $(1,1)$ as in \eqref{1.5}, we chose parameters and initial conditions satisfying hypothesis \eqref{1.3} and \eqref{1.4}.

\begin{figure}[htp]
	\centering
	\includegraphics[scale=0.3]{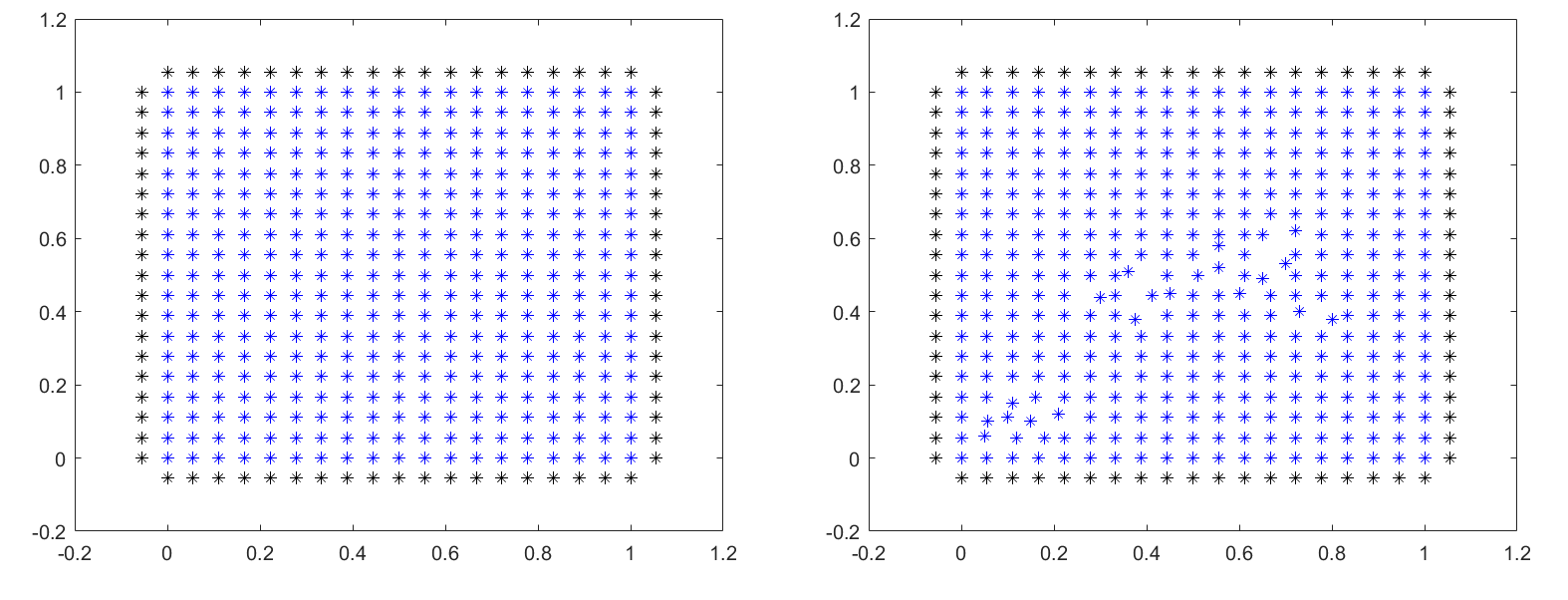}
	\caption{Discretizations of $\Omega$ considered: regular mesh of 19 inner nodes (left) and irregular cloud of points (right).}
	\label{fig:f3}
\end{figure}

For all the examples, we consider the square domain $\Omega = [0,1]\times [0,1]$ with two spatial discretizations, one with a regular mesh and one with an irregular one, to test the differences obtained. Both clouds of points are depicted in Figure \ref{fig:f3}. Notice how along the boundary nodes, colored in black, there is always an inner node placed along the direction of the normal vector. Thus, a simple first order scheme can be used to approximate the normal derivative for the Neumann boundary conditions.

\subsection{Example 1}
\begin{figure}[htp]
    \centering
    \includegraphics[scale=0.3]{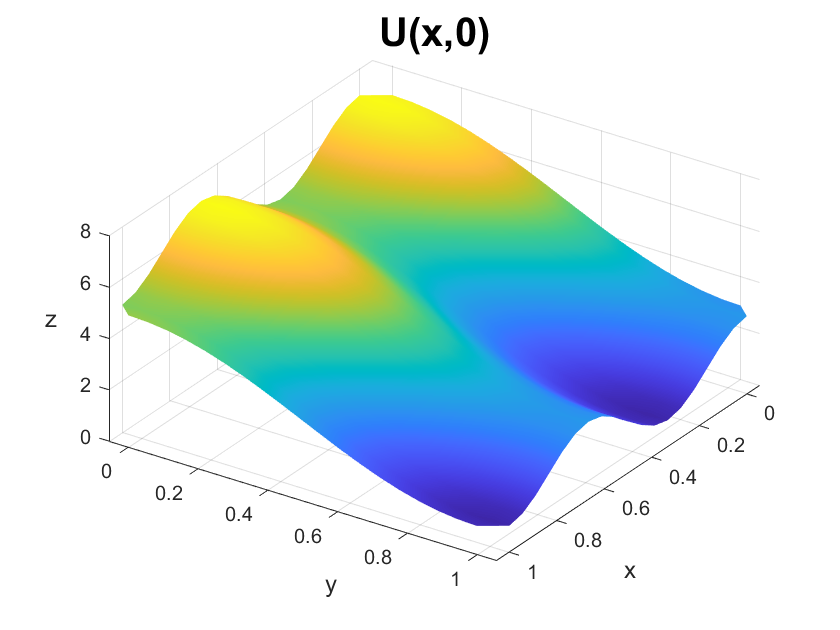}
    \caption{Initial condition for Example 1.}
    \label{fig:f3.2}
\end{figure}

For this first case, as initial condition for $u$ and motility regulation function $\gamma$, we take
\begin{equation} \label{3.1}
    u(\vect{x},0) = 4 + \cos(3 \pi x) + 2 \cos(\pi y), \quad \gamma(s) = e^{-s}.
\end{equation}
Both fulfill conditions \eqref{1.3}-\eqref{1.4} if the logistic growth parameter $\mu$ is taken large enough. In this case, $\mu = 3$ ensures $-2 \gamma'(s) + \gamma''(s) s < \mu$ for all $s \geq 0$. The graph of the initial value $u(\vect{x},0)$ is shown on Figure \ref{fig:f3.2}.
\\
Lastly, as weights for the function $B$ we take $w_{i} := \frac{1}{h_i^2 + k_i^2}$ and a small enough time step, $\Delta t = 0.001$, satisfying the assumption made in Theorem \ref{num1}. 

\begin{figure}[htp]
    \centering
    \includegraphics[scale=0.4]{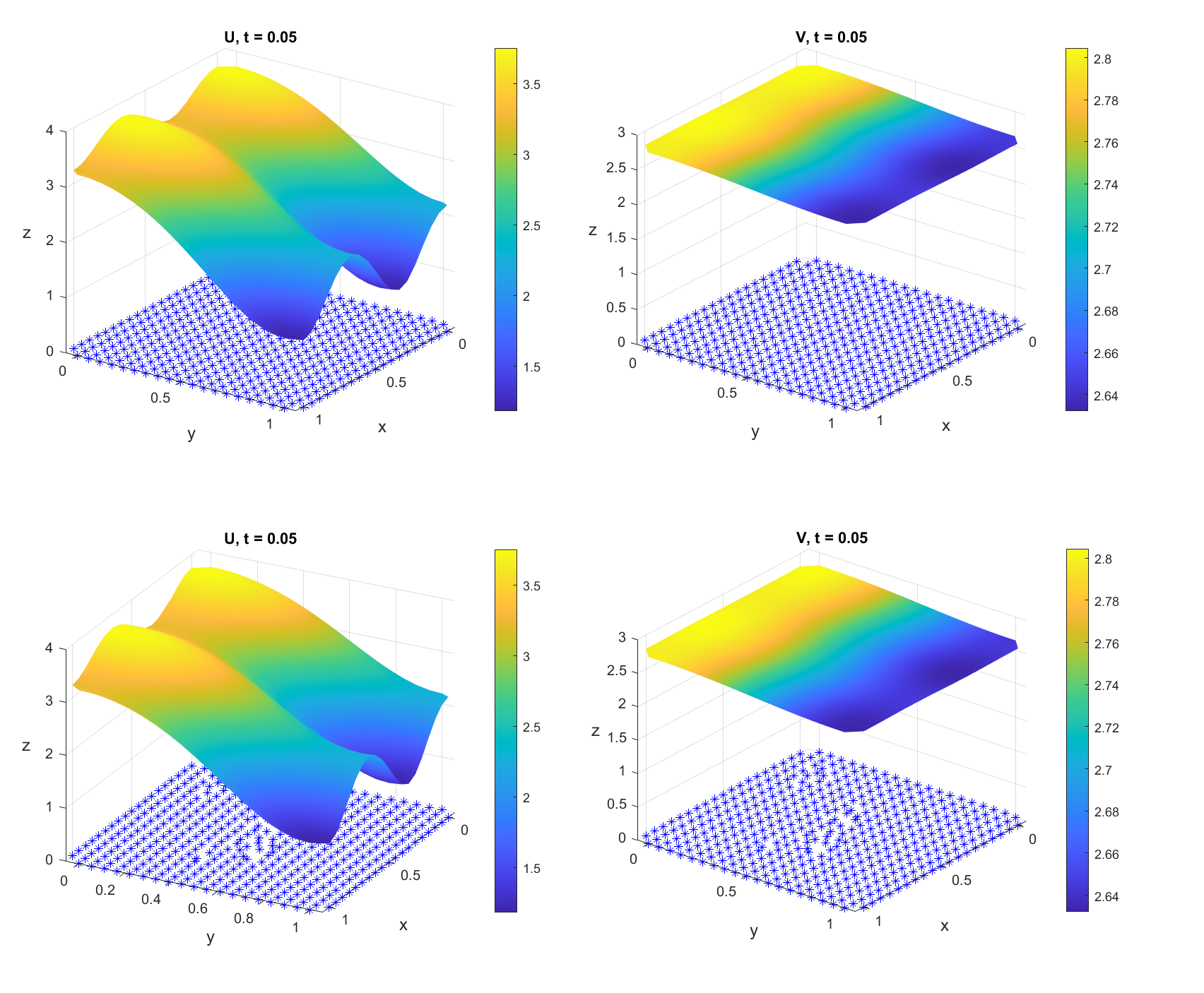}
    \caption{Numerical solution $(U,V)$ at $t = 0.05$ using the regular mesh (upper images) and the irregular cloud of points (lower images).}
    \label{fig:f4}
\end{figure}

The results of the convergence of $||U-1||_{l^\infty(\Omega)}$  and $||V-1||_{l^\infty(\Omega)}$ are gathered in Table \ref{t1}. The convergence is fast, as the growth factor for the logistic model is relatively large.
\begin{table}[h]
\centering
\begin{tabular}{llllll}\hline \\[-1em]
$T(t)$  & 0.05 & 0.1 & 0.5 & 1 & 5 \\ \\[-1em] \hline
\\[-0.5em]
$||U-1||_{l^\infty(\Omega)}$   & 2.7502 & 1.6669 & 0.2086  & 0.0374 & $2.1293 \cdot 10^{-7}$  \\ \\[-0.5em]
$||V-1||_{l^\infty(\Omega)}$    & 1.8045 & 1.2085 & 0.1911  & 0.0368 & $2.1357 \cdot 10^{-7}$  \\ \\[-0.8em] \hline
\end{tabular}
\caption{\label{t1} Values of $||U-1||_{l^\infty(\Omega)}$ and $||V-1||_{l^\infty(\Omega)}$ at different time instants from Example 1, calculated with the regular grid.}
\end{table}
\\
Lastly, the results with both clouds of points are shown in Figure \ref{fig:f4} at time $t = 0.05$, both having very close shape and values.

\subsection{Example 2}
\begin{figure}[htp]
    \centering
    \includegraphics[scale=0.4]{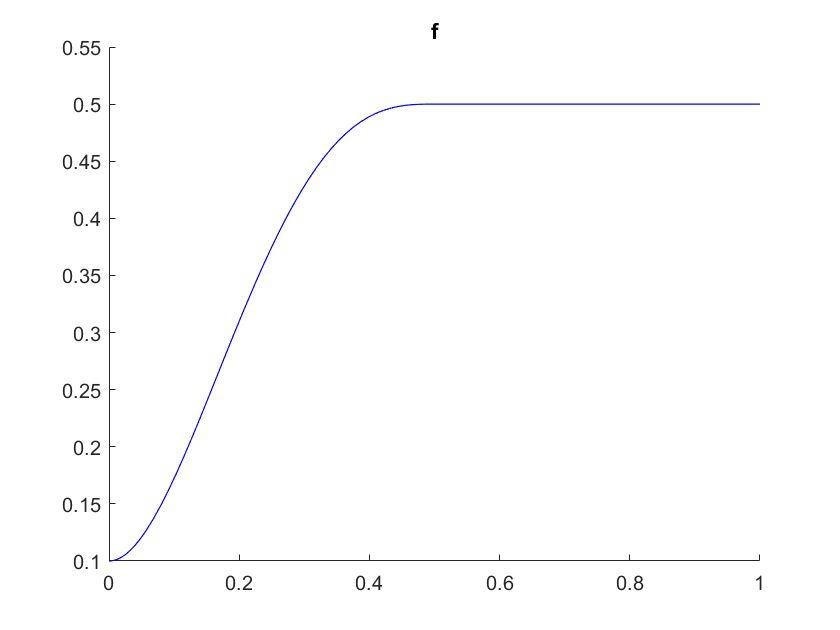}
    \caption{Function $f(x)$ for the initial condition \eqref{3.2}.}
    \label{fig:f5}
\end{figure}

For this second example, we keep the previous domain $\Omega = [0,1] \times [0,1]$, considering the new initial value:
\begin{equation}\label{3.2}
    u(\vect{x},0) = f(x) \cdot [1 + \cos(2 \pi y)],
\end{equation}
where $f$ is a piece-wise function made up of a $4$th degree polynomial for $x \in [0,0.5]$ and a constant function for $x \in (0.5, 1]$, that is:
\begin{equation}\label{3.3}
f(x) = 
    \begin{cases}
        a_4 x^4 + a_3 x^3 + a_2 x^2 + a_1 x + a_0 := p(x), \quad x \in [0,0.5],  \\
        c, \quad x \in (0.5, 1].
    \end{cases}
\end{equation}
To satisfy condition \eqref{1.4} we consider an increasing polynomial with $p(0) > 0 $, $p'(0) = 0$, $p(0.5) = c$, $p'(0.5) = p''(0.5) = 0$. Setting $p(0) = 0.1$ and $c = 0.5$, we have:
$$a_4 = 19.2, \quad a_3 = -25.6, \quad a_2 = 9.6, \quad a_1 = 0, \quad a_0 = 0.1.$$
The graph of $f$ is depicted in Figure \ref{fig:f5}.
 \\
Lastly, we take a new motility regulation function, this time considering $\gamma(s) = (1+s)^{-2}$ and $\mu = 4.5$ to satisfy \eqref{1.3}. The weights and the time step are taken as in Example 1.

\begin{figure}[htp]
    \centering
    \includegraphics[scale = 0.25]{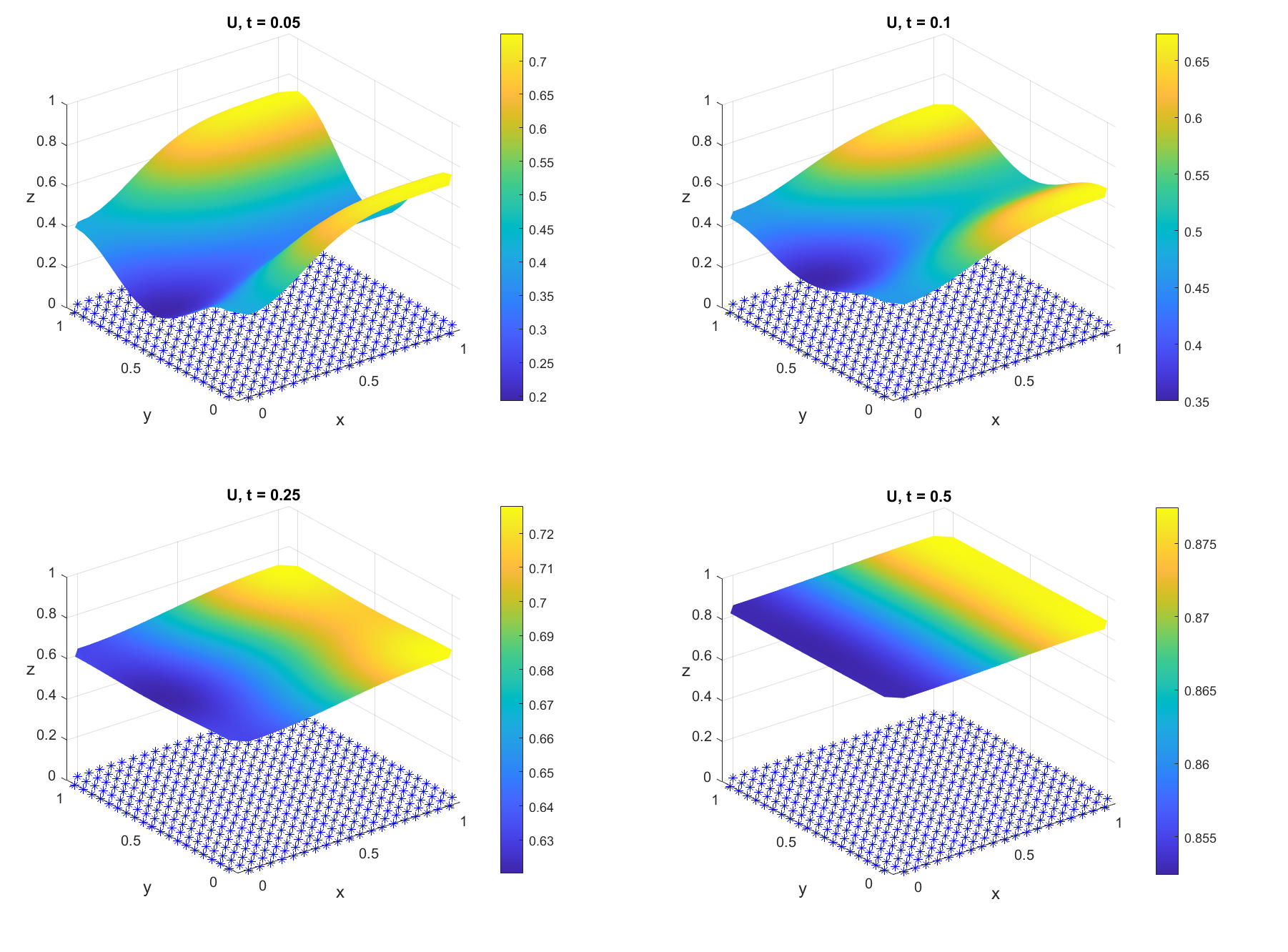}
    \caption{Numerical solution $U$ of Example 2 at different time instants using the regular grid.}
    \label{fig:f6}
\end{figure}

\begin{figure}[htp]
    \centering
    \includegraphics[scale = 0.3]{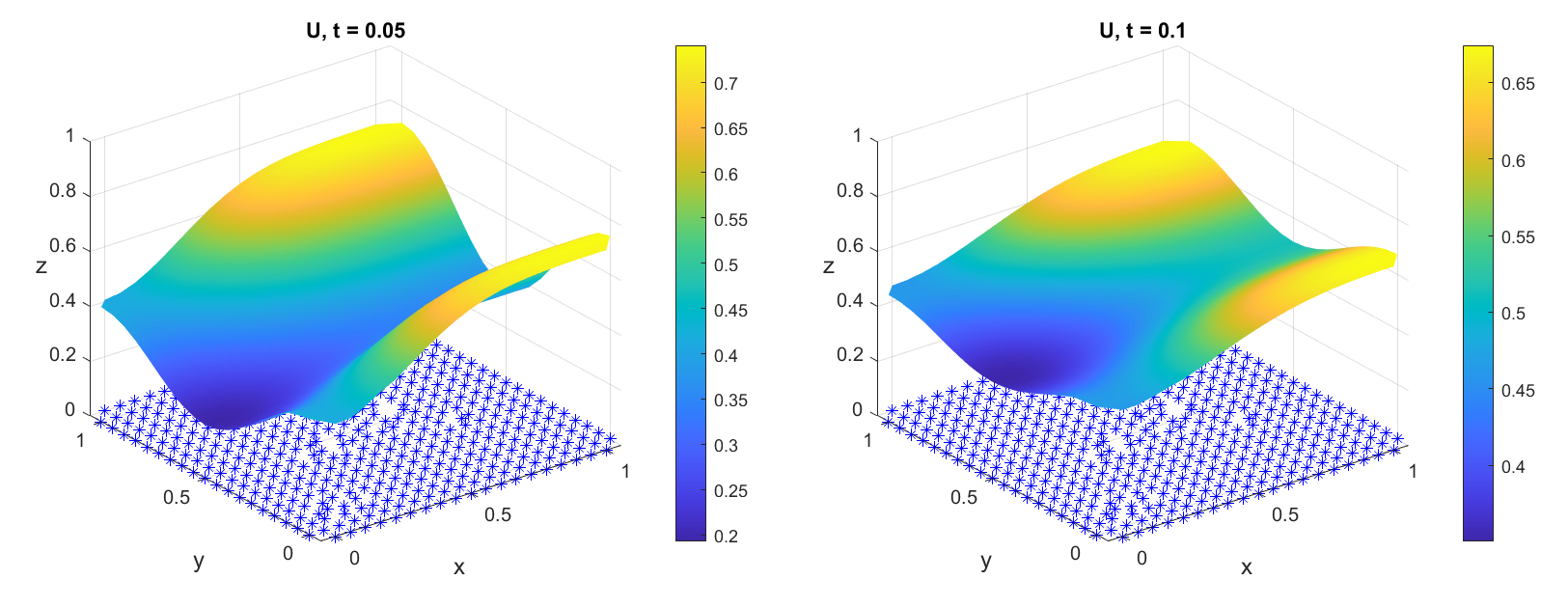}
    \caption{Numerical solution $U$ of Example 2 at different time instants using the irregular grid.}
    \label{fig:f7}
\end{figure}

On Figure \ref{fig:f6} the values of $U$ calculated through the regular grid are displayed at $t = 0.05$, $t = 0.1$, $t = 0.25$ an $t = 0.5$, showing the initial flattening of the solution, afterwards increasing up to $1$, the carrying capacity of the logistic model. 
\\
The results with the irregular cloud of points are nearly identical, as can be seen on Figure \ref{fig:f7}, displaying the values at $t = 0.05$ and $t = 0.1$.
\\
Moreover, the values of $||U-1||_{l^\infty(\Omega)}$  and $||V-1||_{l^\infty(\Omega)}$ are also collected in Table \ref{t2}, verifying the convergence given in \eqref{1.5}.
\begin{table}[h]
\centering
\begin{tabular}{llllll}\hline \\[-1em]
$T(t)$  & 0.05 & 0.1 & 0.5 & 1 & 5 \\ \\[-1em] \hline
\\[-0.5em]
$||U-1||_{l^\infty(\Omega)}$   & 0.8074 & 0.6500 & 0.1476  & 0.0166 & $ 2.3951 \cdot 10^{-10}$  \\ \\[-0.5em]
$||V-1||_{l^\infty(\Omega)}$    & 0.5526 & 0.4950 & 0.1367  & 0.0162 & $2.4264 \cdot 10^{-10}$  \\ \\[-0.8em] \hline
\end{tabular}
\caption{\label{t2} Values of $||U-1||_{l^\infty(\Omega)}$ and $||V-1||_{l^\infty(\Omega)}$ at different time instants from Example 2, calculated with the regular grid.}
\end{table}

\section{Conclusions} \label{s3}
Throughout this paper we've presented an explicit-implicit Generalized Finite Difference Scheme for effectively solving the parabolic-elliptic system \eqref{1.2} of non-linear partial differential equations. After first introducing the basic background of the method in section \ref{preliminares}, the main result was established in Theorem \ref{num1}, obtaining a bound for the time step to guarantee the convergence of the method. Furthermore, an algorithm was provided for the computational implementation of the method.

\begin{figure}[htp]
    \centering
    \includegraphics[scale = 0.5]{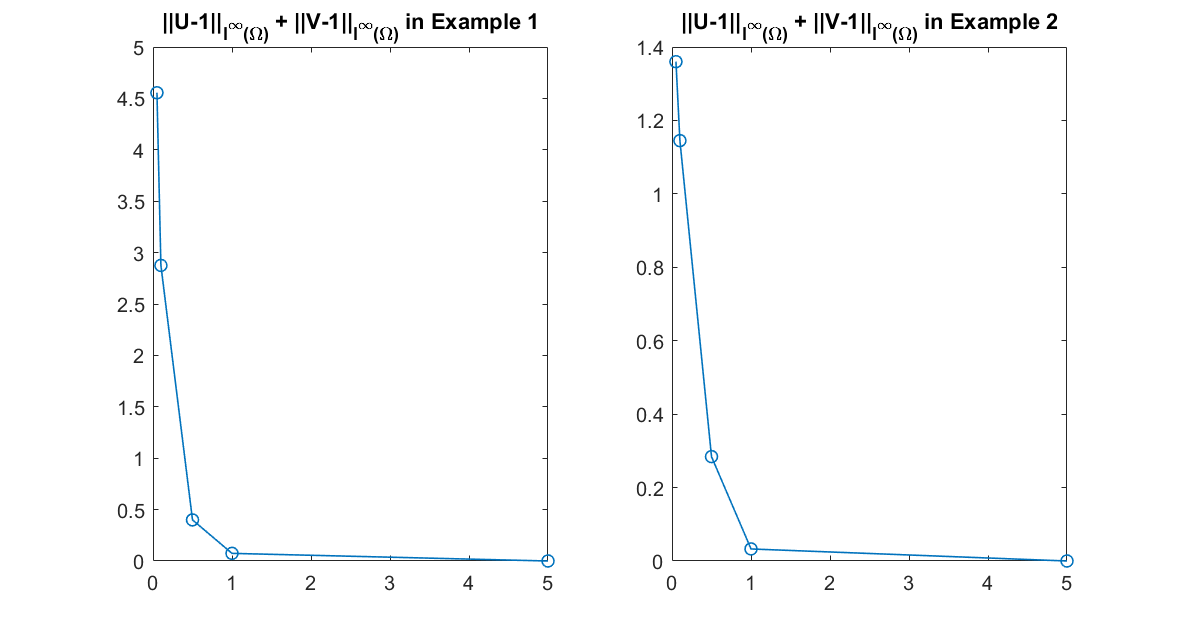}
    \caption{Values of $||U-1||_{l^\infty(\Omega)} + ||V-1||_{l^\infty(\Omega)}$ at different time instants from Example 1 (left) and Example 2 (right).}
    \label{fig:f8}
\end{figure}

Two different examples were tested, over a regular and an irregular cloud of points on the domain $\Omega = [0,1] \times [0,1]$. The results obtained with both discretizations were very close, as seen on Figure \ref{fig:f4} for Example 1 and on Figures \ref{fig:f6} and \ref{fig:f7} for Example 2. The values of $||U-1||_{l^\infty(\Omega)}$ and $||V-1||_{l^\infty(\Omega)}$ were calculated at different time instants verifying the limit \eqref{1.5}. For a graphical representation of the rate of convergence, the values of $||U-1||_{l^\infty(\Omega)} + ||V-1||_{l^\infty(\Omega)}$ for both examples are plotted on Figure \ref{fig:f8}.
\\

\section*{Acknowledgments}
This work was partially supported by the Next Generation European Union funds, the Spanish Ministry of Labour and Social Economy and the Ministry of Economy, Finance and Employment of the Community of Madrid (Grant 17-UCM-INV) (F.H-H.).

\end{document}